\documentclass{amsart} 
\usepackage{graphics} 
\usepackage{fixltx2e}  
\usepackage{amssymb} 
\usepackage{latexsym}
\usepackage{eucal} 
\usepackage{enumerate}
\usepackage{xspace}
\usepackage{newlattice} 

\newtheorem{theorem}{Theorem}

\newcommand{\cornl}[1]{\textup{lc(#1)}}
\newcommand{\cornr}[1]{\textup{rc(#1)}}

\begin{document}
\title[On a result of G\'abor Cz\'edli]
{On a result of G\'abor Cz\'edli\\ concerning congruence lattices\\
of planar semimodular lattices} 
\author{G. Gr\"{a}tzer} 
\address{Department of Mathematics\\ 
   University of Manitoba\\
   Winnipeg, MB R3T 2N2\\
   Canada}
\email[G. Gr\"atzer]{gratzer@me.com}
\urladdr[G. Gr\"atzer]{http://server.math.umanitoba.ca/homepages/gratzer/}

\date{April 28, 1014}

\begin{abstract}
A planar semimodular lattice is \emph{slim} 
if it does not contain $\SM 3$ as a sublattice.
An SPS lattice is a slim, planar, semimodular lattice.

Congruence lattices of SPS lattices satisfy a number of properties.
It was conjectured that these properties characterize them.
A recent result of G\'abor Cz\'edli proves that there is
an eight element (planar) distributive lattice having all these properties that cannot be represented as the congruence lattice
of an SPS lattice.

We provide a new proof.
\end{abstract}

\maketitle 

\section{Introduction}
Let $L$ be a planar semimodular lattice. 
If $L$ is a \emph{slim} lattice (G. Gr\"atzer and E.~Knapp~\cite{GKn08a}), that is, it contains no $\SM 3$ sublattice, 
we call it an SPS lattice (Slim, Planar, Semimodular).

G. Gr\"atzer, H. Lakser, and E.\,T. Schmidt~\cite{GLS98a} prove
that every finite distributive lattice can be represented as the congruence lattice of a planar semidistributive lattice. 
The proof heavily relies on $\SM 3$ sublattices.

I raised  in \cite{gG14a} the problem of characterizing 
congruence lattices of SPS lattices  and
I proved in~\cite{gG14c} the following necessary condition: 

\begin{theorem}\label{T:threecovers}
Let $L$ be an SPS lattice. 
Then the order of join-irreducible congruences of $L$ satisfies the following condition:
\begin{enumerate}
\item[\tup{(CC1)}] an element is covered by at most two elements.
\end{enumerate}
\end{theorem}

G.~Cz\'edli observed that it follows from 
G. Cz\'edli~\cite[Lemma 2.2]{gC13a} that 
the order of join-irreducible congruences of an SPS lattice $L$
also satisfies the following condition:
\begin{enumerate}
\item[\tup{(CC2)}] for any nonmaximal element $a$, 
there are (at least) two distinct maximal elements $m_1$ and $m_2$ above $a$.
\end{enumerate}

G.~Cz\'edli~\cite{gC13a} proved 
(using his Trajectory Coloring Theorem \cite{gC13},
which has a 30 page proof) that the converse does not hold.

\begin{figure}[hbt]
\centerline{\includegraphics{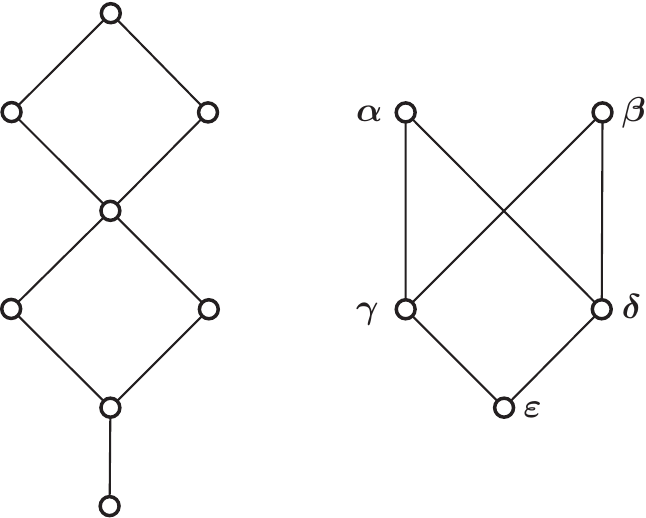}}
\caption{The lattice $D_8$ and the order $P = \Ji(D_8)$}%
\label{Fi:DandP}
\end{figure}

\begin{theorem}\label{T:main}
The eight element distributive lattice $D_8$ of Figure~\ref{Fi:DandP}
cannot be represented as the congruence lattice of an SPS lattice $L$.
\end{theorem}

\begin{figure}[hbt]
\centerline{\includegraphics{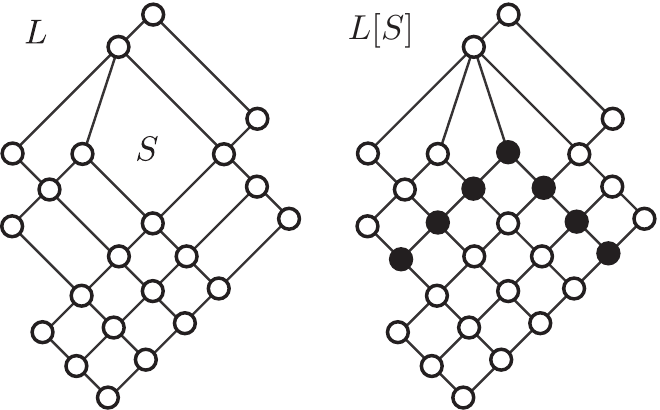}}
\caption{The fork construction, two steps}\label{Fi:forks}
\end{figure}

Note that the order, $\Ji (D_8)$, 
of join-irreducible elements of~$D_8$---see Figure~\ref{Fi:DandP}--- 
satisfies conditions (CC1) and (CC2).

In this note, we provide a new proof of Theorem~\ref{T:main} 
that does not utilize the Trajectory Coloring Theorem. 
However, the new proof utilizes some nontrivial results,
see Theorems~\ref{T:patch} and \ref{T:widetight}.
 
For the basic concepts and notation, see 
G. Gr\"atzer~\cite{LTF}.
For an overview of SPS lattices, see 
G. Cz\'edli and G. Gr\"atzer~\cite{CGa} and
G. Gr\"atzer~\cite{gG13b}, Chapters 3 and 4 of 
G.~Gr\"atzer and F. Wehrung eds. \cite{LTSTA}.

\section{Forks}\label{S:Forks}
We need some basic definitions and results.

Let $L$ be an SPS lattice. A \emph{$4$-cell} in $L$ is a covering 
$\SC 2^2$ with no interior element.

Let $S$ be a $4$-cell of $L$.
We construct a lattice extension $L[S]$ of $L$ as follows. 
 
\begin{figure}[t]
\centerline{\includegraphics{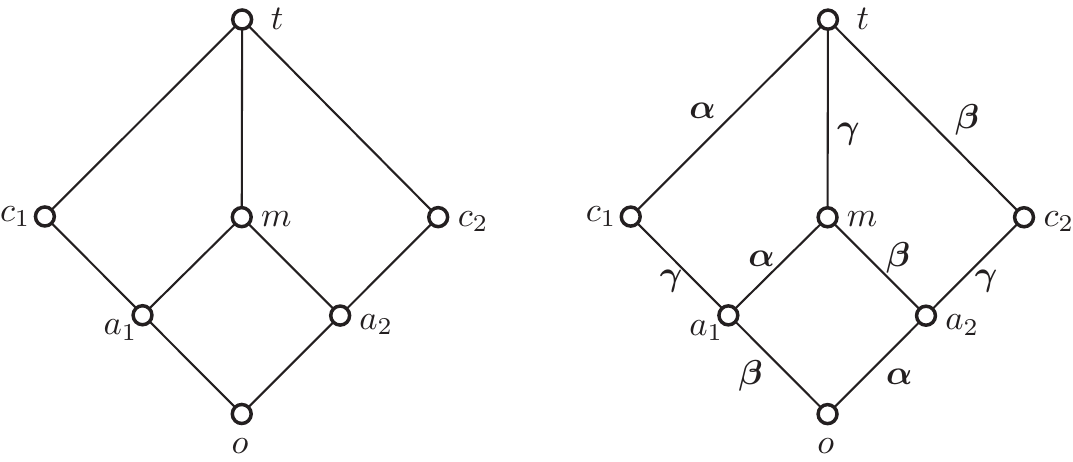}}
\caption{The lattices $\SN 7$ and $L_1$}\label{Fi:n7l1}
\end{figure}

Firstly, we replace $S$ by a copy of $\SN 7$, 
the lattice of Figure~\ref{Fi:n7l1}, introducing three new element. 

Secondly, we do a series 
of steps---each step introducing one new element: 
if there is a chain $u \prec v \prec w$ such that~$v$ 
is a new element but $u$ and $w$ are not,
and  $T = \set{u \mm x, x, u, w = x \jj u}$ 
is a $4$-cell in the original lattice~$L$, 
see Figure~\ref{Fi:forks}, then we insert a~new element~$y$ 
such that $x \mm u \prec y \prec x \jj u$ and $y \prec v$. 
Figure~\ref{Fi:forks} shows some steps of the construction.

Let $L[S]$ denote the lattice we obtain when the procedure terminates.
We say that $L[S]$ is obtained from
$L$ by \emph{inserting a fork} at the $4$-cell~$S$.

Let $L$ be a planar lattice. A \emph{left corner} (resp., \emph{right corner}) of $L$ is a doubly irreducible element in $L - \set{0,1}$ on the left (resp., right) boundary of~$L$. 
As~in G.~Gr\"atzer and E.~Knapp~\cite{GKn09}, 
we define a \emph{rectangular lattice} $L$ 
as a planar semimodular lattice 
which has exactly one left corner, $\cornl L$, 
and exactly one right corner, $\cornr L$, 
and they are complementary.
A rectangular lattice $L$ is a \emph{patch lattice} 
(G. Cz\'edli and E.\,T. Schmidt~\cite{CSa} and \cite{CSb}) if
$\cornl L$ and $\cornr L$ are dual atoms.

We need the following result of G. Cz\'edli and E.\,T. Schmidt~\cite{CSa}:

\begin{theorem}[Structure Theorem for SPS Lattices]
\label{T:patch}
An SPS lattice $L$ can be
obtained from a planar distributive lattice $D$ by inserting forks.
\end{theorem}

We call the planar distributive lattice $D$ the \emph{grid}
of the SPS lattice $L$. 

In an SPS lattice $L$, we call the covering square of 
$S =\set{o, c_l, c_r, t}$  a~\emph{tight square}, 
if $t$ covers exactly two elements, namely, $c_l$ and $c_r$, in $L$;
otherwise, $S$ is a \emph{wide square}.

The following results are from G. Gr\"atzer~\cite{gG14a} (the notation
$o$, $a_l$, $a_r$, $t$, and $m$ is from Figure~\ref{Fi:n7l1}, denoting the elements of the inserted $\SN 7$).

\begin{theorem}\label{T:widetight}
Let $L$ be an SPS lattice and let $S$ be a covering square of $L$. 
If $S$ is a
\begin{enumeratei}
\item \emph{wide square}, 
then the congruence $\bgg(S) = \consub{L[S]}{m,t}$ of $L[S]$ is generated by a~congruence of $L$.
\item \emph{tight square}, 
then $L[S]$ has exactly one join-irreducible congruence, 
namely $\bgg(S) = \consub{L[S]}{m,t}$,
that is \emph{not generated} by a~congruence of $L$. 
\end{enumeratei}

\end{theorem}

\section{Proof of Theorem~\ref{T:main}}\label{S:proof}
Let $L$ be a finite SPS lattice whose join-irreducible congruences
form an order~$P$ as in Figure~\ref{Fi:DandP}. 
By G. Gr\"atzer and E. Knapp~\cite{GKn09}, 
$L$ has a congruence-preserving extension
to a rectangular lattice, 
so we can assume that $L$ is a rectangular lattice.

\emph{Case 1:} $L$ is a patch lattice.

By the Structure Theorem for SPS Lattices, we can obtain 
$L$ from $\SC 2^2$ by inserting forks. 
Let $\SC 2^2 = L_0, L_1, \dots, L_n = L$ 
be a sequence of fork insertions from $\SC 2^2$ to $L$. 

There is only one way to insert
a fork into $\SC 2^2$; 
so $L_1$ is the lattice of Figure~\ref{Fi:n7l1}.
There is one more join-irreducible congruence in $L_1$, 
the congruence~$\bgg$. 
We ``color'' the diagram of $L_1$
with the join-irreducible congruences generated by the edges.

\begin{figure}[t!]
\centerline{\includegraphics{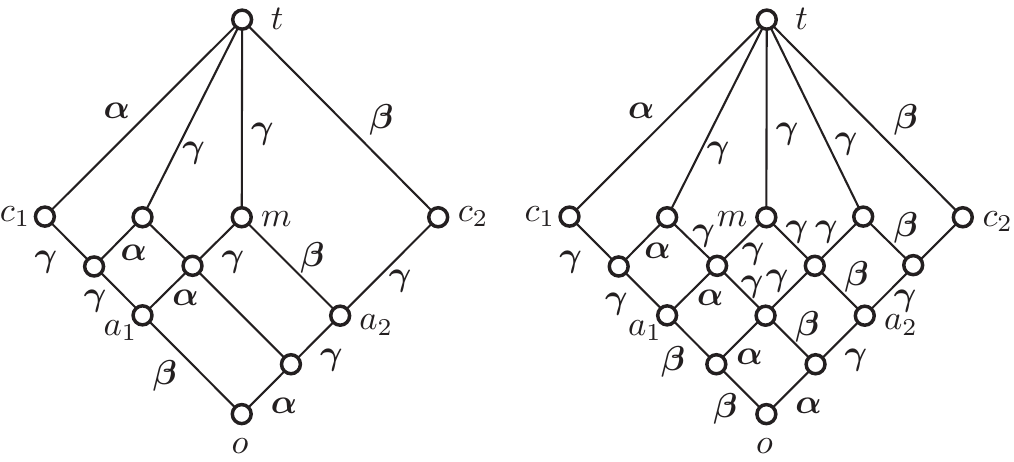}}
\caption{The lattices $L_2$ and $L_3$ with no new congruence}\label{Fi:L2L3nonewcong}
\end{figure}

To get $L_2$, we pick a covering square $S_1$ in $L_1$ and
insert a fork into $L_1$ at $S_1$. 
If~the top element of $S_1$ is $t$, then $S_1$ is a wide square,
so we get no new congruence by Theorem~\ref{T:widetight}(i), 
see Figure~\ref{Fi:L2L3nonewcong}. The next step, again with no new congruence is shown in the same figure.
Note that all covering squares 
of the lattices $L_1$, $L_2$, $L_3$ satisfy the following condition:
\begin{enumeratei}
\item[(Col)] all covering squares are colored by $\bgg$ by itself
or with $\bga$ or with $\bgb$ with one exception: 
the bottom covering square is colored by $\bga$ and $\bgb$.
\end{enumeratei}

We proceed thus in $k-1$ steps to get $L_k$, where $k \geq 1$ is the largest number with the property that the number of join-irreducible 
congruences do not change from $L_{k-1}$ to $L_k$. 
Note that $k < n$ because $L_n = L$ has more join-irreducible congruences. 
Clearly, $L_k$ also satisfies condition (Col).

We proceed to $L_{k+1}$ by inserting a fork 
into a covering square $S$ of $L_k$.
There are five cases:

\begin{enumeratei}

\item the top element of $S$ is $t$;

\item the top element of $S$ is not $t$ 
and $S$ is ``monochromatic'', colored by $\bgg$;

\item the top element of $S$ is not $t$ 
and $S$ is colored by $\set{\bgg,\bga}$;

\item the top element of $S$ is not $t$ and
$S$ is colored by $\set{\bgg,\bgb}$;

\item $S$ is the bottom covering square 
colored by $\set{\bga,\bgb}$. 
\end{enumeratei} 

Case (i) cannot happen, it would contradict Theorem~\ref{T:widetight}(i)
and the definition of~$k$.

If Case (ii) holds, $S$ is tight, 
so by Theorem~\ref{T:widetight}(ii),
we add a join-irreducible congruence $\bgg' < \bgg$. 
By~Figure~\ref{Fi:DandP}, we must have $\bgg' = \bge$.
This is a contradiction because no fork insertion can add an element 
$\bgd$ between two existing elements.

Case (iii) proceeds the same was as Case (ii).

Case (iv) is symmetric to Case (iii).

\begin{figure}[t!]
\centerline{\includegraphics{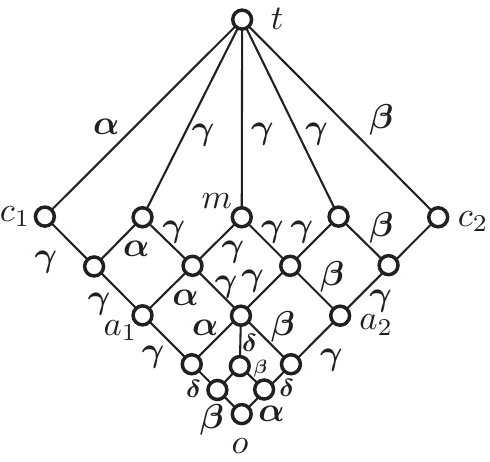}}
\caption{The lattice $L_{k+1}$ for case (v)}\label{Fi:Lk}

\bigskip

\bigskip

\centerline{\includegraphics{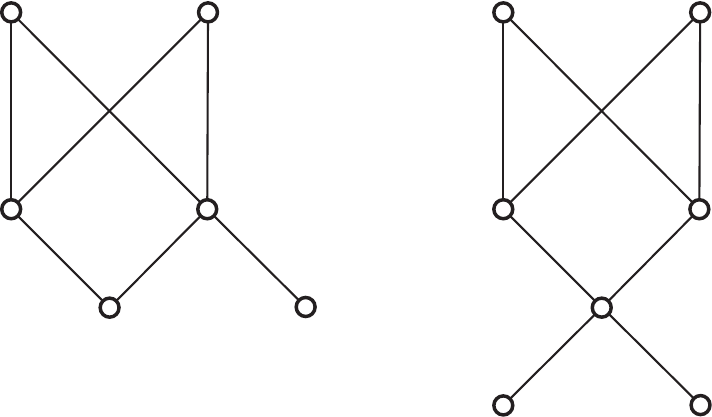}}
\caption{Larger nonrepresentable orders}\label{Fi:largenonrep}
\end{figure}

So we are left with Case (v). 
In this case, the top element of $S$ is not $t$, so $S$ is tight.
By Theorem~\ref{T:widetight}(ii), 
we get a new join-irreducible congruence $\bgd$ in $L_{k+1}$  
satisfying that $\bgd < \bga$ and $\bgd < \bgb$, see Figure~\ref{Fi:Lk}. 

Let $k+1 \leq m < n$ be the largest integer so that 
the number of join-irreducible congruences
does not change from $L_{k+1}$ to $L_m$---we insert forks into wide squares. 
Then the lattices $L_{k+1},\dots, L_{m}$ share 
the following property of $L_{k+1}$: 
\begin{enumeratei}
\item[(P)] there are only  ``monochromatic'' squares and 
$\set{\bga, \bgb}$, $\set{\bga, \bgg}$, 
$\set{\bga, \bgd}$, $\set{\bgb, \bgg}$, $\set{\bgb, \bgd}$ squares, 
but there is no $\set{\bgg, \bgd}$ square.
\end{enumeratei}

We obtain the last join-irreducible congruence, $\bge$, in $L_{m+1}$. 
By property (P), we cannot have $\bge \prec \bgg$ and $\bge\prec \bgd$.

\emph{Case 2:} $L$ is not a patch lattice.
It follows that the grid $D$ is not a covering square.
The lattice $D$ has three join-irreducible congruences 
forming an antichain of $3$ elements. 
But the order of Figure~\ref{Fi:DandP} has only 
one- and two-element antichains.

This completes the proof of Theorem~\ref{T:main}.

The two cases of the proof can easily unified. 
If $L$ is not a patch lattice, 
then there are a number of covering squares we can start with.
Say, we start with the covering square $S$ colored by $\bga$ and $\bgb$. 
We proceed as in Case 1, except that there are more colors. 
The additional colours do not change the fact that 
there are no $\set{\bgg, \bgd}$ squares.

\section{Larger nonrepresentable orders}\label{S:Discussion}
The proof of Theorem~\ref{T:main} 
can be made to work under related assumptions. 
For instance, the orders in Figure~\ref{Fi:largenonrep}
are not representable.

\end{document}